Boris R. Vainberg (on his 80th birthday)

Boris R. Vainberg was born on March 17, 1938, in Moscow. His father was a Lead Engineer in an aviation design institute. His mother was a homemaker. From early age, Boris was attracted to mathematics and spent much of his time at home and in school working through collections of practice problems for the Moscow Mathematical Olympiad. His first mathematical library consisted of the books he received as one of the prize-winners of these olympiads.

In 1955 Boris became a student and then a graduate student at the Differential Equations Division of the Mechanics and Mathematics Department of Moscow State University. His PhD adviser was S.A.Galpern. The head of the Differential Equations Division at that time was Ivan Georgievich Petrovsky, who was also the rector of Moscow State University. Petrovsky can be viewed as a symbol of what we now call "The Golden Age of Mechanics and Mathematics Department". The faculty of the Differential Equations Division at that time included V.I.Arnold, M.I.Vishik, E.M.Landis, O.A.Oleinik. Other divisions of the department also had some of the greatest mathematicians of the 20th century: A.N.Kolmogorov (Probability Theory), P.S.Aleksandrov (Topology), I.M.Gelfand and D.E.Menshov (Theory of Functions and Functional Analysis), this list can go on and on.

In those years of political liberalization, the Mechanics and Mathematics Department and the Division of the Differential Equations were at the forefront of the newest scientific ideas. Those included the theory of generalized functions started by S.L.Sobolev in the 1930-s, extended by L. Schwartz in the 1950-s, and further developed by I.M.Gelfand, E.G.Shilov, and others. The most modern results of the time included the breakthrough work by L.Hormander on pseudo-differential operators, integral Fourier operators, and hypoellipticity. Some of the early fundamental papers by B.R.Vainberg (join with V.V.Grushin) concerned this very subject: the analysis of Dirichlet-to-Neumann map from the point of view of pseudo-differential operators (see below).

In 1963 Boris defended the PhD thesis and, with essential support of I.G.Petrovsky, was offered an assistant professor position at the Department. B.R.Vainberg worked in the Division of Differential Equations for 28 years. His contribution was essential in the creation of the Correspondence School at the Mechanics and Mathematics Department, organized by I.M.Gelfand. Boris attracted students from the Mechanics and Mathematics Department to grade the assignments of the students from the Correspondence School, and all the interaction between them was conducted under his supervision. In August of 1991 Boris Vainberg moved to the US. During the first year he was a visiting professor at the University of Delaware, and since 1992 has been a professor at the University of North Carolina at Charlotte.

Boris Vainberg is an author of three monographs [1-3], a chapter in another book [4], and around 170 papers, mostly in the areas of mathematical physics and partial differential equations.

The results of his PhD thesis "Conditions at infinity providing for the unique solvability of hypoelliptic equations in the entire space" [5] provide the radiation conditions of Sommerfeld type for general elliptic operators (in particular, for elasticity equations). Much later, these results were applied in his paper with W. Shaban (his PhD student at the time) to study the discrete Laplace equation [6], where radiation conditions have many peculiarities: there exist several scattered waves, and the limiting absorption principle is violated for some frequencies.

In 1968, jointly with V.V. Grushin, he became a prize-winner of the Moscow Mathematical Society (the annual prize for young mathematicians) for the work on uniformly noncoercive problems for elliptic equations [7, 8]. One of the important applications of this work was the result that the Dirichlet-to-Newman map is a pseudo-differential operator, with the evaluation of its complete symbol.

In 1970 B. R. Vainberg presented the doctoral (advanced PhD) dissertation "Elliptic problems in exterior domains and large time behavior of the solutions of the Cauchy problem for hyperbolic equations". The dissertation was rejected due the anti-semitic policies that were rapidly gaining strength since the late 1960-s. When the dissertation of an excellent mathematician G.I.Eskin (currently a UCLA professor) was rejected a little earlier, many viewed this as a random event. After the rejection of Vainberg's dissertation, it was realized that the situation in the department took a drastic turn for the worse.

In 1973 B. R. Vainberg, together with V. Maz'ya, described the steady-state oscillations of a fluid layer of variable depth caused by the periodic or uniform motion of a body immersed in the fluid [9, 10]. They found the geometrical conditions on the inhomogeneity (elevation of the bottom of the water layer and the shape of the body) that provide the absence of the eigenvalues imbedded into the continuous spectrum and lead to the unique solvability of the problem. For operators with periodic coefficients, the absence of the embedded eigenvalues was established by Bosis later in a joint work with P. A.Kuchment [11]. In 1981, B.R.Vainberg, jointly with V.G.Mazya, studied the characteristic Cauchy problem for general hyperbolic equations [12]. Ten years later, their result was re-discovered by L.Hormander, but only for equations of the second order.

In 1987 B. Vainberg presented another doctoral thesis "Local Energy Decay for Exterior Hyperbolic Problems and Quisi-Classical Approximations", and was awarded the degree of doctor of physical and mathematical sciences. In the thesis, he developed a direct method [11, 12, 1, 2] to obtain the large time asymptotic behavior of the local energy and of the solutions to non-stationary problems in the exteriors of non-trapping domains. The method is based on high frequency estimates and low frequency asymptotics of the solutions of the corresponding stationary problems. In particular, these results provide all the consequences that follow from the Lax-Phillips scattering theory. Later he extended his approach to the case of time periodic media and obstacles [15-16]. These results of B. Vainberg on the decay of local energy were successfully used by him and the co-authors in the study of asymptotic stability of stationary states in nonlinear wave equations, and in the study of solitons for a particle interacting with its wave field and the Klein-Gordon field [17-18].

After moving to the U.S., B Vainberg obtained many profound results in collaboration with S. Molchanov. They found spectral asymptotics for operators in domains with fractal boundaries [19-20], with sparse potentials [21-23], and for other important classes of operators [24-26], studied the spectrum of the Schrodinger operator with slowly decaying and random potentials [27-33]. A series of their papers [34-37] concerns the propagation of waves in complex networks of thin fibers, with applications to fiber optics. They found the asymptotics of solutions when the thickness of fibers tends to zero, recovered the gluing conditions at vertices of the one-dimensional limiting graph, and used this simplified problem on the graph to describe the propagation of waves in the original problem. Together with co-authors, they studied mathematical models for homopolymers [38,39] and the spectral properties of non-local Schrodinger operators [40], established the global limit theorems for random walks with "heavy tails", and applied them to describe the front propagation in biological models and to study intermittency in these models [41].

B. Vainberg wrote a series of papers [42-44] with E.Lakshtanov on interior transmission eigenvalues (an object that appears in the scattering by an obstacle). In particular, they obtained a new Weyl law, where the eigenvalues are counted with the plus or minus signs that are defined by the direction of motion of the corresponding eigenvalues of the scattering matrix. In another work (joint with R. Novikov), they showed that the global Riemann-Hilbert problem can be applied to solve two-dimensional inverse scattering problems in all the cases, in particular, in the case when exceptional points are present [46]. This allowed them to solve certain important nonlinear equations of soliton theory in dimension 2+1 (such as the focusing Davey-Stewartson equation) without the usual assumption on the smallness of initial data [45-48]. B. Vainberg made important contributions to many other problems of mathematical physics.

Boris has a wonderful family, he keeps playing tennis and work on mathematical problems. Let us wish him excellent health and new successes in mathematics.

Yu. Egorov, A. Komech, P. Kuchment, E. Lakshtanov, V. Mazya, S. Molchanov, R. Novikov, M. Freidlin.

Cited papers by B.R.Vainberg.

1. Asymptotic Methods in Equations of Mathematical Physics, 1982, Moscow, Moscow Univ. Publishers (in Russian).

2. Asymptotic Methods in Equations of Mathematical Physics (revised and expanded English version), Gordon and Breach Science Publishers, New York--London, 1989.

3. Linear Water Waves: A Mathematical Approach, Cambridge University Press, 2002 (with N. Kuznetsov and V. Maz'ya).

4. Large Time Asymptotic Expansion of the Solutions of Exterior Boundary Value Problems for Hyperbolic Equations and Quasiclassical Approximations, Chapter in "Partial Differential Equations, V", Series: Encyclopaedia of Math. Sciences, Springer-Verlag, Berlin-Heidelberg-New York, 1999.

5. Principles of radiation, limiting absorption and limiting amplitude in the general theory of partial differential equations, Russian Math. Surveys, Vol 21, No 3, pp 115-193, (1966).

6. Radiation conditions for the difference Schrodinger operators, Applicable Analysis, 80, 525-556 (2001) (with W. Shaban).

7. Uniformly nonelliptic problems I, Math. USSR Sbornik, Vol 1, No 4, pp 543-568, (1967), (with V. Grushin).

8. Uniformly nonelliptic problems II, Math. USSR Sbornik, Vol 2, No 1, pp 111-133, (1967), (with V. Grushin).

9. On the problem of the steady--state oscillations of a fluid layer of variable depth, Trans. Moscow Math. Soc., Vol 28, pp 56--73, (1973), (with V. Maz'ya).

10. On the plane problem of the motion of a body immersed in a fluid, Trans. Moscow Math. Soc., Vol 28, pp 33-55, (1973), (with V. Maz'ya).

11. On absence of embedded eigenvalues for Schrodinger operators with perturbed periodic potentials, Comm. PDE, 25, 1809-1826, (2000) (with P. Kuchment).

12, The characteristic Cauchy problem for a hyperbolic equation. (Russian. English summary) Trudy Sem. Petrovsk. No. 7 (1981), 101–117), (with V. Maz'ya).